\documentclass{amsart}
\usepackage{amsmath,amsthm,amscd}

\newtheorem{thm}{Theorem}[section]
\newtheorem{cor}[thm]{Corollary}
\newtheorem{lem}[thm]{Lemma}
\newtheorem{prop}[thm]{Proposition}
\newtheorem{conj}[thm]{Conjecture}

\newtheorem{question}{Question}

\theoremstyle{definition}

\newtheorem{notation}{Notation}
\newtheorem{rem}{Remark}[section]

\newcommand{\cL}{\mathcal{L}}

\newcommand{\cO}{\mathcal{O}}
\newcommand{\cR}{\mathcal{R}}

\newcommand{\fa}{\mathfrak{a}}
\newcommand{\fb}{\mathfrak{b}}

\newcommand{\fm}{\mathfrak{m}}

\newcommand{\fM}{\mathfrak{M}}
\newcommand{\fN}{\mathfrak{N}}

\newcommand{\bbI}{\mathbb{I}}
\newcommand{\bbN}{\mathbb{N}}
\newcommand{\NN}{\mathbb{N}}
\newcommand{\VV}{\mathbb{V}}
\newcommand{\bbZ}{\mathbb{Z}}
\newcommand{\bbQ}{\mathbb{Q}}

\newcommand{\bbP}{\mathbb{P}}

\newcommand{\gh}{\eta}
\newcommand{\gl}{\lambda}

\newcommand{\lrarrow}{\longrightarrow}
\newcommand{\la}{\longrightarrow}

\DeclareMathOperator{\SSpec}{\mathcal S {\mathit p  \mathit e \mathit c}}

\DeclareMathOperator{\Hom}{Hom}

\DeclareMathOperator{\Ann}{Ann}

\DeclareMathOperator{\Sym}{Sym}
\DeclareMathOperator{\core}{core}
\DeclareMathOperator{\grcore}{gradedcore}

\DeclareMathOperator{\adj}{adj}

\DeclareMathOperator{\Ker}{Ker}
\DeclareMathOperator{\Coker}{Coker}

\DeclareMathOperator{\Proj}{Proj}
\DeclareMathOperator{\Spec}{Spec}

\begin{document}

\title[]{Core versus graded core, and global sections of line bundles}
\author{Eero Hyry}
\address{Department of Mathematics,
University of Helsinki\\
 Helsinki, Finland}
\email{Eero.Hyry@helsinki.fi}
\thanks{Research of Hyry supported by the National Academy of 
Finland, project number 48556}
\author{Karen E. Smith}
\address{Department of Mathematics, University of Michigan, Ann Arbor, MI,
48109-1109}
\email{kesmith@umich.edu}
\thanks{Research of Smith partially supported by the Clay Foundation and by 
the US National Science Foundation Grant DMS 00-70722.}


\maketitle

\section{Introduction}

When does an ample line bundle on a smooth projective variety have a
non-zero global section?
In this paper, we show that this question is equivalent to
 a fundamental problem in commutative algebra regarding the equality of
 {\it core}  and {\it graded core}  for 
 a certain associated  homogeneous ideal. This builds on the
 project begun in \cite{HS}, where it was shown that a sufficiently good
algebraic understanding of the graded core can be used to show the existence
of global sections of line bundles. This paper essentially 
treats the converse: we show  the existence of  
non-zero sections gives rise to nice formulas for cores.

The core of an ideal in a
Noetherian  commutative 
 ring  is
the intersection of all its reductions---that is, the
intersection of all subideals having the same integral closure.
 The core first arose in the work of Rees and Sally \cite{RS}
because of its connection with Brian\c con-Skoda theorems,
and has recently been the subject of active investigation in 
commutative algebra; see \cite{HunSwan, CPU, CPU2, HS, PU, HT}.

For 
 a homogeneous ideal in a graded ring,
it is also natural to consider its {\it graded core},
namely, the intersection of all its {\it homogeneous} reductions.
 The core and graded core of a homogeneous ideal are 
both homogeneous ideals and there is an obvious inclusion 
 of the core in the graded core.
 A natural 
question is: when does equality hold? This question arose in 
the work \cite{HS} in finding sections of line bundles, and has also been 
considered by Huneke and Trung in \cite{HT} for purely algebraic reasons.

Quite generally, the core and graded core are equal for
homogeneous ideals  generated by  elements of the same degree
(see Lemma \ref{standardlem} for a  precise statement).
On the other hand, there are examples of ideals
having so few homogeneous reductions that it is easy to see the core is 
strictly smaller. 
 But  what can be said for ideals having many 
graded reductions?

In this paper we answer this question for a specific type of ideal in 
a {\it section ring}, while providing one answer to the opening question
 about sections of line bundles. 
Specifically, we prove a formula for the graded core
 (see Theorem \ref{gradform} and Corollary \ref{cor1})
 which, when combined
with the formula for core from \cite{HS}, yields the following result.

\begin{thm} 
\label{main}
Let $\cL$ be an ample invertible sheaf on an irreducible 
projective variety $X$ over a field of characteristic zero,
and let \begin{equation*}S=\bigoplus_{n\ge 0}H^0(X, \cL^n)\end{equation*} be 
the corresponding 
section ring. Assume that $S$ is Gorenstein.
Fix $N\gg 0$, and let $I=S_{\ge N}$ be the ideal 
generated by all elements of degree at least $N$.
Then  $\grcore(I)=\core(I)$ if and only if $\cL$ has a non-trivial
 global section.
Furthermore, in this case,  
$$\grcore(I)=\core(I) = S_{\ge Nd + a + 1}$$
where $d$ is the dimension of $S$ and $a$ is its $a$-invariant.{\footnote{By definition, $a$ invariant of $S$ is defined by 
$a = - \min \{n\, | \, [\omega_S]_n \neq 0$\}.}}
 \end{thm}

This theorem makes it easy to construct examples of 
 homogeneous ideals with abundant
 homogeneous reductions where  core and graded core
differ. 
On the other hand, 
 one might also  hope to use  the theorem to algebraically prove the 
existence of global  sections of particular line bundles---
indeed, this is our 
motivation   for studying core. 

Theorem \ref{main} follows from a more general result which
is valid without the Gorenstein assumption on $S$.
 To motivate this, recall that in \cite{HS} it was shown that 
  from the point of view of projective geometry, 
 one should really look at
cores inside the
canonical module of $S$. That is, the object in which  we are interested
is the module $\grcore(J\omega_S)$,
defined as the intersection in $\omega_S$ of the submodules $J\omega_S$
as $J$ ranges over all homogeneous reductions of $I$.
With this definition, we prove the 
 following, which specializes immediately to  Theorem \ref{main}
if $S$ happens to be Gorenstein.

\begin{thm} 
\label{mainthm}
Let $\cL$ be an ample invertible sheaf on an irreducible Cohen-Macaulay
projective variety $X$ of positive
dimension $d-1$  over a field of characteristic zero,  
and let \begin{equation*}S=\bigoplus_{n\ge 0}H^0(X, \cL^n)\end{equation*} be 
the corresponding 
section ring. 
Fix $N\gg 0$, and let $I=S_{\ge N}$ be the ideal 
generated by all elements of degree at least $N$.
Then $H^0(X,\cL)\not=0$ if and only if $\grcore(I\omega_S)=\core(I\omega_S)$.
Furthermore, in this case, we have
$$
\grcore(I\omega_S)=\core(I\omega_S) = [\omega_S]_{\geq dn+1}.
$$
 \end{thm}

\bigskip
Our methods also yield the following formula for core in standard graded rings
of {\it arbitrary} characteristic.

\begin{thm}\label{standard}
Let $(S, \fm)$ be a standard graded reduced Cohen-Macaulay
 ring of dimension $d$ over an infinite
 field of arbitrary characteristic, and let 
$a$ denote its $a$-invariant.
Then for all  $N \geq 1 $, 
$$
\core(\fm^N) = \grcore(\fm^N) =  \fm^{Nd + a + 1}.
$$
\end{thm}

Similar formulas appear in \cite[\S 6]{HS}, 
 \cite{PU} and \cite{HT}. Note that  Theorem
\ref{standard} does not require
 the ground field to have characteristic zero, as in these other formulas.
Indeed, according to \cite[Example 5.9]{PU}, in general,
the core of the maximal ideal can depend on the characteristic.

\medskip
Our original motivation for studying core was the following 
remarkable conjecture of   Kawamata{\footnote{This was also
 posed as a question by Ambro in \cite{Am}.}:
 an ample line bundle
$\cL$ on a smooth variety $X$ always has a non-zero global section provided
that $\cL \otimes \omega_X^{-1}$ is also ample.
In \cite{HS}, we showed that this conjecture  follows
 from a conjectured formula for graded core generalizing a formula due to 
Huneke and Swanson.
As a corollary to the main theorem of the current paper, 
we see that in fact this conjectured formula for the graded core is
{\it  equivalent}  to Kawamata's Conjecture; 
see Conjecture \ref{KC1}.

\smallskip
The contents of the paper are as follows.
Section 2 contains a list of properties about section rings.
The main work is done in 
 Section 3, where the main result,  Theorem \ref{gradform}, is proven.
Section 4
contains the proof of
  Theorem \ref{standard}.
Section 5 contains
a list of commutative algebra questions and detailed remarks motivating them.
 Progress
towards any of these questions
may lead to progress 
 towards Kawamata's conjecture or a number of other 
related conjectures in birational algebraic geometry.
 We believe
 that commutative algebraists ought to be studying section rings
more, and we hope our
 questions help stimulate research in this direction.
We also hope that the geometric connections
will help other researchers to
appreciate the core.

\section{Section Rings}

Let $X$ be a projective scheme over a field
 and let $\cL$ be an ample
invertible sheaf on $X$. The {\it section ring}
of the pair $(X, \cL)$ is the $\bbN$-graded ring
$$
S = \bigoplus_{n \in \bbN} H^0(X, \cL^n).
$$
This  is a finitely generated graded ring, and
there is a natural isomorphism $\Proj S \cong X$ under which  the
sheaves on $\Proj S$ given by the graded module $S(n)$ correspond to
the invertible sheaves $\cL^n$. (Here $S(n)$ denotes the graded cyclic 
$S$-module generated in degree $-n$.) The ring $S$ is reduced if and only if 
$X$ is reduced, and is a domain if and only
if $X$ is reduced and irreducible. In this case, furthermore,
$S$ is normal if and only if $X$ is normal. 
Moreover, 
 the sheaf 
given by the graded canonical module $\omega_S$ of $S$ corresponds to 
$\omega_X$.

Section rings are much nicer than arbitrary graded rings.
They share many of the properties of standard graded rings (that is,
those generated in degree one) and should be
thought of as a natural generalization of this standardly considered case.
A normal section ring is standard graded if and only if the line bundle 
$\cL$ is very ample and projectively normal;
thus the passage from standard graded to section rings in commutative
algebra  is 
 analogous to
the passage from very ample to ample line bundles in algebraic geometry.
Under this analogy, the powers of 
the homogeneous maximal ideal of a standard graded ring 
correspond to the ideals
generated by elements of degrees at least $N$ in a section ring.
It is these ideals whose cores are crucial in understanding sections of $\cL$.

For future reference, we summarize here a few of the nice properties of
section rings. All of these properties can fail for an arbitrary 
graded ring.

\begin{prop}[Properties of Section Rings]
Let $S$ be the section ring of a projective variety,
and let $\fm$ denote its unique homogeneous maximal ideal.
Then 
\begin{enumerate}
\item The section ring $S$ has depth at least two (or depth one in the
case of a section ring of dimension one);  
that is $H^0_{\fm}(S) = 0$  and, if $\dim S > 1$,   $H^1_{\fm}(S) = 0 $ 
as well. 
\item
 The sheaf corresponding to the graded module
 $S(1)$ is an invertible sheaf on $\Proj S$.
\item 
The sheaves corresponding to the graded modules
 $S(n)$ are invertible on  $\Proj S$ for all $n$.
\item
There exist homogeneous non-zero-divisors $x_1, \dots, x_r$ generating an 
$\fm$-primary ideal such that the $\bbZ$-graded ring 
$S_{x_i}$ is isomorphic 
to the ring $[S_{x_i}]_0[t, t^{-1}]$, 
where $t$ is an indeterminate of degree one and $[S_{x_i}]_0$
is the subring of the localization ring $S_{x_i}$
consisting of degree zero elements.
\item For any homogeneous element $x$ and any $n\in \bbZ$,
the $n$-th graded piece of $S_x$ is a rank one projective
module over the subring $[S_x]_0$.
\item For all $n \gg 0$, the elements of degree $n$ 
generate an $\fm$-primary ideal.
\item For all $n \gg 0$, $S_n^p = S_{np}$ for all $p\geq 1$.
\item For all $p, q \gg 0$, $ S_p S_q = S_{p+q}.$
\item The ring $S$ has isolated singularity (respectively, 
isolated non-rationally singular point, 
isolated non-Cohen-Macaulay point, isolated
non-Gorenstein point, et cetera)
if and only if the variety $\Proj S$
 is non-singular (respectively, rationally singular, Cohen-Macaulay, 
Gorenstein, et cetera.)
\item The Hilbert function $n \mapsto \dim S_n$ is polynomial for $n \gg 0$.
\item If $\Proj S$ is smooth (respectively, has rational singularities, 
is Cohen-Macaulay, et cetera), then also $\Proj S^{\natural}$ is smooth
(respectively, has rational singularities, 
is Cohen-Macaulay, et cetera), where $S^{\natural} = 
\bigoplus_{n \geq 0}S_{\geq n}$ is the Rees ring
of the natural filtration $\left\{S_{\geq n}\right\}_{n \in \bbN}$.
\end{enumerate}
Furthermore, 
Properties 
2,  3, and 4   each characterize normal
section rings
among all finitely generated normal graded domains.
\end{prop}

\begin{proof}
We only sketch the proofs since all can be found in the
literature.
The original source for section rings is \cite[Section 4.5]{EGAII}.
In \cite{Demazure}, normal graded rings are studied quite generally
and shown to be rings of sections for
 ample $\bbQ$-Cartier Weil divisors---
 the section ring case corresponds to the divisor actually being Cartier. 

Property 1 follows by computing $H^i_{\fm}(S)$ from the 
extended \v Cech complex
 $$
0 \la S \la \oplus S_{x_i} \la 
\oplus_{i < j} S_{x_ix_j} \la \dots
$$
where $x_0, \dots, x_d$ is a homogeneous system of parameters for $S$.
This leads to an exact sequence
$$
0 \la 
H^0_{\fm}(S) \la S \la \oplus_{n \in \bbZ} H^0(X, \cL^n) \la H^1_{\fm}(S)
\la 0.
$$
 By definition of $S$, then, we see that both $H^0_{\fm}(S)$ and 
$H^1_{\fm}(S)$ must vanish. (The computation when $S$
 has dimension one  is similar).

Properties 2 and 3 are immediate since the sheaf associated to 
$S(n)$ is isomorphic to $\cL^n$, which is invertible. These properties
characterize  section rings by the main theorem of \cite{Demazure}.

Property 4 is essentially a restatement of 3. The $x_i$ are chosen so that 
the open sets $U_i = \Spec [S_{x_i}]_0$ in $\Proj S$ where $x_i$ does not
vanish are a trivializing cover  for $\cL$. Then $\cL(U_i) = 
[S_{x_i}]_n$ is a free  $\cO_{\Proj S}(U_i) = [S_{x_i}]_0$-module of
rank one. Letting $t$ be a generator for $\cL(U_i)$, it follows
that $t^n$ is a generator for $\cL^n$. The isomorphism follows.
Conversely, if the isomorphism holds, then each $\cL^n(U_i)$ is clearly
free of rank one.  
Property 5 is essentially Property 4 stated for any arbitrary cover of
$\Proj S$ rather than only for a cover that trivializes $\cL$.
See also \cite{Sm1}.

Properties 6, 7 and 8  follow from the definition of ampleness.
In particular, Property 6 follows because $\cL^n$ is globally generated
for all large $n$: then $\cL^n$ has a set of global sections that 
do not simultaneously vanish, which means that the ideal of $S$ they generate
is $\fm$-primary.
 Property 7 follows since $\cL^n$ is very ample and projectively normal
for all large $n$. Property 8 follows from 
the K\"unneth formula; see, for example, \cite[Example 1.2.18]{PAG}.

Property 9 follows immediately from Property 4 because the map
$A \lrarrow A[t, t^{-1}]$ is smooth for any ring $A$.

Property 10 follows from Serre Vanishing.
Indeed, the function $n \mapsto  \chi(X, \cL^n)$ is polynomial
in $n$ (its precise form being given by the Riemann-Roch formula),
where $\chi(X, \cL^n) =
 h^0(X, \cL^n) - h^1(X, \cL^n) + \dots \pm h^d(X, \cL^n)$
\cite{Sn}\cite{Kleiman}.
For $n \gg 0$, Serre vanishing ensures that all the higher cohomology
terms are zero. 

Finally, Property 11 is essentially treated, for example,
 in \cite[Paragraph 6.2.1]{HS} or \cite[Section 8.7.3]{EGAII}.
 The point is that $\Proj S^{\natural}$
is isomorphic to the total space of the line bundle $\cL$ over $X$.
Thus there is a smooth map $\Proj S^{\natural} \lrarrow \Proj S$  making
all the nice properties of $\Proj S$ pass to $\Proj S^{\natural}$.
\end{proof}

We will use most of these properties without explicit mention in the
next section.  In fact, the only properties we
 really use are (6), (7) and (8),
but because section rings are a natural and important class, we prefer to
state our results for them.

\section{The Main Theorem}

For a  ring $S$, recall that $\core(I\omega_S)$ denotes the
intersection of all submodules $J\omega_S$ of $\omega_S$ where $J$ is a
 reduction of $I$. Similarly, when $I$ is a homogeneous ideal in a graded  ring
$S$, 
$\grcore(I\omega_S)$ is the corresponding  intersection
 over all homogeneous
reductions $J$ of $I$.
Our goal is to prove the following theorem.

\begin{thm} 
\label{gradform}
Let $S$ be an equidimensional section ring of dimension $d \geq 2$ and
characteristic zero, and assume that $\Proj S$ is Cohen-Macaulay. 
 Fix $N\gg 0$, and let $I=S_{\ge N}$ be the ideal 
generated by all elements of degree at least $N$.
{\footnote{The precise value required for $N$ is not very important;
 what is used is that $\cL^n$ is very ample,
 and  that the higher cohomology modules
 of $\omega_X\otimes \cL^n$ vanish for all 
$n \geq N$.}}
Then, 
$$
\grcore(I\omega_S) = \bigoplus_{i \in \bbZ; S_i = 0}
[\omega_S]_{Nd-i}
$$
as a graded submodule of $\omega_S$.
 \end{thm}

Note that in particular (taking negative values of $i$),
 $\grcore(I\omega_S)$ contains
the submodule $[\omega_S]_{\geq Nd+1}$ generated by elements of 
degrees at least $Nd+1$. By \cite[Corollary 6.4.1]{HS},
this submodule is precisely $\core(I\omega_S)$. But 
if $S$ is ``missing'' a component of degree $i$, then the graded
core picks up elements of degree $Nd-i$. In particular, Theorem \ref{mainthm}
from the introduction is an immediate corollary of Theorem \ref{gradform}.

\bigskip
 We now turn to the proof of Theorem \ref{gradform}.
 We begin with a lemma that reduces this problem to 
a core computation for a slightly different ideal.

\begin{lem}\label{redI'} 
Let $S$ be an equidimensional
 section ring such that $\Proj S$ is Cohen-Macaulay.
  Fix $N\gg 0$, and 
let $I$ denote the ideal generated by elements of degree
at least $N$. Then $\grcore (I\omega_S) =\core (I'\omega_S)$ where 
$I'$ denotes the ideal of $S$ generated by all elements of degree precisely
$N$. 
\end{lem}

\begin{proof} Because $N\gg 0$, we know that $S$ has systems
 of parameters of degree $N$. In
particular, $I'$ is a reduction of $I$, and we must have
 $\grcore (I\omega_S) \subset \grcore (I'\omega_S)$. 
On the other hand, every 
homogeneous minimal reduction $J$ of 
$ I$ is generated by elements of degree $N$ 
(see \cite[Proposition 2.1.3]{HS}).
 So $\grcore (I\omega_S) =\grcore (I'\omega_S)$. 
But as $I'$ is generated by 
homogeneous elements 
of the same degree $N$,  Lemma~\ref{standardlem} below ensures 
that  $\core (I'\omega_S) =\fa_1\omega_S\cap \ldots \cap \fa_r\omega_S$ 
where $\fa_1,\ldots,\fa_r$ are homogeneous 
minimal reductions of $I$. To apply Lemma~\ref{standardlem}, note that
since both modules are $\fm$-primary, we can check equality after localization.
The hypothesis  that $S_{\fm}$ is a generalized
Cohen-Macaulay ring is satisfied because $S$ is equidimensional with (at worst)
an
isolated non-Cohen-Macaulay point.
 The proof that 
  $\grcore (I'\omega_S)=\core (I'\omega_S)$ is complete.
\end{proof}

\begin{lem} 
\label{standardlem}
Let $(A,\fm)$ be a generalized Cohen-Macaulay local ring of 
dimension $d$ containing its residue field $k=A/\fm$. 
Let $I\subset A$ be an $\fm$-primary ideal. Let 
$\{a_1,\ldots,a_{\mu}\}$ be a set of generators for $I$. 
If $I\subset \fm^N$ where $N\gg 0$, then there exist minimal reductions 
$\fa_1,\ldots,\fa_t$ of  $I$ such that $\core(I\omega_A)
=\fa_1\omega_A\cap \ldots \cap \fa_t\omega_A$ and 
for each $\fa_h$ there is a generating set $\{b_{h1},\ldots, b_{hd}\}$ where 
\begin{equation*}b_{hi}=\sum_{j=0}^{\mu} \gl_{hij}a_j\
 \quad {\text{ with }} \quad \gl_{hij}\in k \quad {\text{for}}
\quad h=1,\ldots,t;\ 
i=1,\ldots,d.
\end{equation*}
A similar statement holds also for cores of ideals rather than
cores inside the canonical module.
\end{lem}

\begin{proof} This is a 
slight generalization of a
 special case of ~\cite[Theorem 4.5]{CPU}. We must remove the
hypothesis there that $A$ is Cohen-Macaulay and make some minor
adjustments  
to replace the core ideal by the core submodule in
$\omega_A$.

Because there exists an $r$ such that 
 every reduction of $I$ contains $I^{r+1}$,
it follows that   $\omega_A/\core(I\omega_S)$ has finite length. So there
 are in any case finitely many 
reductions $J_1,\ldots,J_t$ of $ I$ such that $\core(I\omega_A)=
J_1\omega_A\cap \ldots \cap J_t\omega_A$. The Cohen-Macaulayness 
of $A$ was used in 
 the proof of ~\cite[Theorem 4.5]{CPU} only  
  in the proof of~\cite[Lemma 4.3]{CPU}, which 
 established
the equality $l_A(I\omega_A/\fa_i\omega_A)=l_A(I\omega_A/J_i\omega_A)$.
We will instead establish this equality under the assumption that $A$ is
a generalized Cohen-Macaulay ring.

Let $M$ be a generalized Cohen-Macaulay $A$-module,
meaning that the local cohomology modules $H^i_{m}(M)$ are all
of finite length, for $i < \dim M$.
Recall (see, for example, \cite{T}) that an $\fm$-primary ideal
$\fb$ of $ A$ is called standard with respect to $M$ if 
every system of parameters $(b_1,\ldots,b_d)$  for $A/\Ann M$ 
contained in $\fb$
satisfies
$l_A(M/(b_1,\dots,b_d)M)-e(b_1,\ldots,b_d;M)=\bbI(M)$, where $\bbI(M)$
 is the Buchsbaum invariant $\sum_{i=0}^{\dim M -1}\binom{\dim M - 1}{i}
l_A(H^i_{\fm}(M))$
of 
$M$. By~\cite[Satz 3.3]{CST},  there always exists an
 $N\gg 0$ such that every ideal 
$\fb$ contained in $\fm^N$
 is standard with respect to M.
 Our assumption 
that $N$ is large, together with~\cite[Theorem 3.17]{GY}, 
imply that $I$ is standard with respect to $\omega_A$. So 
\begin{align*}
\begin{split}
l_A(I\omega_A/J_i\omega_A)&=l_A(\omega_A/J_i\omega_A)-l_A(\omega_A/I\omega_A)\\
&=e(J_i;\omega_A)+\bbI(\omega_A)-l_A(\omega_A/I\omega_A)\\
&=e(\fa_i;\omega_A)+
\bbI(\omega_A)-l_A(\omega_A/I\omega_A)=l_A(I\omega_A/\fa_i\omega_A)
\end{split}
\end{align*}
as needed. 
\end{proof}

\bigskip
To compute 
 $\core(I'\omega_S)$, we will
apply the main technical theorem of 
\cite{HS}.  First, we need to recall some notation from Section 3.1 of that
paper.

\medskip
\begin{notation}[Adjoint Type Modules]\label{adj}
Let $(A,\fm)$ be a local
 ring which is a homomorphic image of a Gorenstein local ring, and let
  $I$ be any   ideal of $A$ of positive height.
 Set $Y=\Proj A[It]$ where $A[It]$ is the Rees algebra of
$I$.
 For each integer $p$,  we define the  $p$-th
{\it  adjoint type module }
by 
\begin{equation}\label{p}
\Omega_p= \Gamma(Y,I^p\omega_Y).
\end{equation}
Note that the adjoint type modules depend on the ideal $I$,
although this dependence is suppressed from the notation.

Clearly, 
 $\Omega_{p+1}\subset \Omega_p$
 and $I\Omega_p \subset \Omega_{p+1}$ 
for all $p\in \bbZ$. 
Recall that $\Omega_0$ can be considered as a submodule of $\omega_A$ 
by means of the trace homomorphism $\Omega_0 \lrarrow \omega_A$ for the
blowing up map $Y \rightarrow \Spec A$. 
(When $I$ has height at least two,
the trace homomorphism can be defined in an elementary way as follows:
Let $U\subset Y$ be the non-empty open set of $Y$ isomorphic under the
blowing up to $\Spec A \setminus \VV(I)$.
Then 
$$
\Omega_0 = \Gamma(Y,\, \omega_Y) \overset{restr.}\la  \Gamma(U,\, \omega_Y)
\cong \Gamma(\Spec A \setminus \VV(I),\, \omega_A) = \omega_A,$$
with the first arrow being the restriction map and the last equality 
holding because $\omega_A$ satisfies Serre's $S_2$ condition.)
In particular, we
 get an $I$-filtration $\omega_A \supset \Omega_1\supset \Omega_2 \dots$\ .
As explained  in \cite[Paragraph 2.6.2]{HS}, for positive values of $p$,
 there is a natural identification of
$\Omega_p$ with the $p$-th graded component of the graded canonical 
module  $\omega_{A[It]}$ of $A[It]$. 
\end{notation}

\bigskip

The proof of Theorem~\ref{gradform} is based on the following result:

\begin{thm}
\label{mainalg}
Let $S$ be a finitely generated $\bbN$-graded domain
of dimension $d \geq 2$ containing the rational numbers, and  let $\fm$ denote
its homogeneous maximal ideal.
 Let $I$ be an $\fm$-primary ideal generated by elements of all 
the same degree. Assume that for any homogeneous reduction $J$ 
of $I$,  
\begin{equation}
\label{eq1}
J\omega_S \cap \Omega_{d-1}
=J\Omega_{d-2}
\end{equation}
 as submodules of $\omega_S$,
where $\Omega_p$ is the $p$-th adjoint type module
of $I$ as defined in (\ref{p}).   Then $\grcore(I\omega_S)\subset \Omega_d$. 
\end{thm} 

\begin{proof}
This is a graded version of the main technical theorem 
of \cite{HS} and is
direct consequence of Lemma 3.3.1 and Theorem 3.4.1 there.
\end{proof}

Thus combining Lemma \ref{redI'} and Theorem \ref{mainalg}
we see that to  compute the graded core of $I\omega_S$, 
we need to  verify (\ref{eq1}) for the adjoint type modules associated to 
$I'$,  where  $I'$ is  
the ideal of $S$ generated by the elements of degree $N$.
We use the notation
 $$\Omega'_p=\Gamma(Y',{I'}^p\Omega_{Y'})$$
 where $Y'=\Proj S[I't]$ for these modules.

\medskip

\begin{notation}[$S^{\natural}$ and $S^{\dagger}$]\label{rings} 
Recall the ``natural" construction from~\cite[Section 8.7.3]{EGAII}. 
Given a graded ring $S$, we form
the graded ring 
\begin{equation*}S^{\natural}=\bigoplus_{p\ge 0}S_{\ge p}\end{equation*} 
where
 $S_{\ge p}$ is the ideal generated by   all elements
of 
$S$ of degrees at least $p$. 
The ring $S^{\natural}$ is the Rees ring corresponding to the
``natural" filtration $\left\{S_{\ge p}\right\}_{p\ge 0}$,
and the corresponding associated graded ring  recovers $S$.
Note that if $I=S_{\ge N}$ where $N\gg 0$ is chosen 
in such a way that $S_{\geq pN}=(S_{\geq N})^p$ for all $p\ge 1$, 
then the Veronese subring $(S^{\natural})^{(N)}$
 coincides with the Rees algebra $S[It]$. 

We define a graded subring $S^{\dagger}$ of $S^{\natural}$ by setting
 \begin{equation*}
S^{\dagger}=\bigoplus_{p\ge 0}S_pS
\end{equation*}
 where $S_pS$
 denotes the ideal of $S$ generated by elements of degree precisely $p$.
 If $N\gg 0$ is chosen as above, then the Veronese
subring  $(S^{\dagger})^{(N)}$ coincides with the 
Rees ring $S[I't]$, where $I'$ is the ideal of $S$ generated by the elements
of degree exactly $N$.

Note that $S^{\natural}$ is finite and birational over 
 $S^{\dagger}$ 
(and if $S$ is normal, then $S^{\natural}$ is the normalization of
$S^{\dagger}$).
Both rings 
have dimension $d +1$, where $d$ is the dimension of $S$.
Furthermore, both rings have  natural 
 bigradings:
\begin{equation*}
S^{\dagger}=\bigoplus_{p,q\ge 0}S_pS_q\qquad\text{and}\qquad
S^{\natural}=\bigoplus_{p,q\ge 0}S_{p+q}.
\end{equation*}
We will make much use of these bigradings. There is also a third grading,
namely, the
 ``total degree,''
whose $n$ graded components are defined as follows:
\begin{equation*}
[S^{\dagger}]_n=\bigoplus_{p,q\ge 0, p+q = n}S_pS_q\qquad\text{and}\qquad
[S^{\natural}]_n=\bigoplus_{p,q\ge 0, p+q = n}S_{p+q}.
\end{equation*}
\end{notation}

We  compute the $p$-th adjoint type module
$\Omega_p'$ of $I'$  by interpreting it as the
appropriate graded piece of the graded canonical module of 
$S^{\dagger}.$  With this in mind, the point is 
prove  the following theorem.

\begin{thm} 
\label{canonicalmodulethm}
Let $S$ be an equidimensional section ring of dimension $d\ge 2$, with
$\Proj S$ Cohen-Macaulay.
Let $G$ denote the set  of  (positive and negative) integers $n$ for which 
$S_n$ is zero. 
Then  for any positive $p$
such that  $[H^{d-1}_{\fm}(S)]_n = 0 $ for all $n < -p$,
we have 
\begin{equation*}
[\omega_{S^{\dagger}}]_p= 
\bigoplus_{-q\in G}[\omega_S]_{p+q}
\end{equation*}
Furthermore, if $S$ has negative  a-invariant,
then 
we also get the following description for non-positive $p$:
\begin{equation*}
[\omega_{S^{\dagger}}]_p= 
\begin{cases}
{[\omega_S]_{\ge p+1}}&\text{if $-p\in G$;}\\
0&\text{otherwise.}
\end{cases}
\end{equation*}
\end{thm}

\begin{rem}
Note that the set $G$ of Theorem \ref{canonicalmodulethm}
 contains all the negative integers. In particular, 
$[\omega_{S^{\dagger}}]_p$ contains the submodule ${[\omega_S]}_{\geq p+1}$.
Since $G$ contains only finitely many positive integers, 
$[\omega_{S^{\dagger}}]_p$ differs from  ${[\omega_S]}_{\geq p+1}$
only by  a vector space of finite dimension.
\end{rem}

\begin{rem}
If $S$ is Cohen-Macaulay or if $\Proj S$ has rational singularities
(and is over a field of characteristic zero), then 
 $[H^{d-1}_{\fm}(S)]$ vanishes  in all negative degrees,
so the cohomological hypothesis of Theorem \ref{canonicalmodulethm}
   is satisfied automatically.
(When $\Proj S$ has rational singularities, this
  follows from the  Kodaira Vanishing Theorem; 
see \cite{HunSm}.) 
Furthermore, when $\Proj S$ is Cohen-Macaulay, it is satisfied for 
all sufficiently large $p$.
\end{rem}

\begin{proof} 
We will study the graded pieces of $\omega_{S^{\dagger}}$ by comparison
with those of $\omega_{S^{\natural}}$, which are easier
to understand.

First note that the finite birational map 
$S^{\dagger} \hookrightarrow S^{\natural}$ induces a natural inclusion
(the trace map) 
$\omega_{S^{\natural}} \hookrightarrow \omega_{S^{\dagger}}$.
We wish to understand the cokernel. We do this by studying the
kernel of the dual map.

Let $\fM$ and $\fN$ denote the homogeneous maximal ideals 
of $S^{\natural}$ and $S^{\dagger}$ respectively.
Note that 
the radical of $\fN S^{\natural}$ is $\fM$.
 Since $\dim S^{\natural}/S^{\dagger}
<d+1$, we obtain a surjective homomorphism\begin{equation*}
\psi \colon H_{\fN}^{d+1}(S^{\dagger})
 \lrarrow H_{\fN}^{d+1}(S^{\natural})=H_{\fM}^{d+1}(S^{\natural}).
\end{equation*}
This surjection is dual to the trace map 
$\omega_{S^{\natural}}\hookrightarrow \omega_{S^{\dagger}}$.
Thus we need to understand the 
 kernel of $\psi$. 

To this purpose, set \begin{equation*}\fN'=\bigoplus_{p\ge 1,q\ge 0}S^{\dagger}_{p,q}, \quad
\fN''=\bigoplus_{p\ge 0,q\ge 1}S^{\dagger}_{p,q}
\quad\text{and}\quad \fN^+=\bigoplus_{p\ge 1,q\ge 1}S^{\dagger}_{p,q}
.\end{equation*}  Observe that $\fN'+\fN''=\fN$ and $\fN'\cap \fN''
=\fN^+$. The 
Mayer-Vietoris sequences for the cohomology
 of $S^{\dagger}$ and $S^{\natural}$ 
respectively  
yield the 
commutative diagram
\begin{equation*}
\begin{array}{ccccccccccc}
& &H_{\fN'}^d(S^{\dagger})\oplus H_{\fN''}^d(S^{\dagger})&\lrarrow&H_{\fN^+}^d(S^{\dagger})&\lrarrow&H_{\fN}^{d+1}(S^{\dagger})&\lrarrow&0\\
& &\varphi\Big\downarrow                         &            &\Big\downarrow & &\Big\downarrow\psi             & & \\
H_{\fN}^d(S^{\natural})&\lrarrow&H_{\fN'}^d(S^{\natural})\oplus H_{\fN''}^d(S^{\natural})&\lrarrow&H_{\fN^+}^d(S^{\natural})&\lrarrow&H_{\fN}^{d+1}(S^{\natural}).& &\\
\end{array}
\end{equation*}

Now,
the module $H^d_{\fN}(S^{\natural})$ vanishes in negative
$p$-degrees (and hence in negative $q$-degrees by symmetry).
 Indeed,  the Sancho de Salas sequence (see \cite{L2} or
 \cite[Paragraph 2.6]{HS})
yields an exact sequence of graded modules
$$
\dots \bigoplus_{p\in \bbZ}H^{d-1}_E(\cO_Y(p)) \lrarrow 
H^d_{\fN}(S^{\natural}) \lrarrow \bigoplus_{p \in \bbN} H^d_m(S_{\geq p}) \lrarrow
\dots,
$$
where $Y = \Proj S^{\natural}$ and  $E$ is the closed fiber of 
the blowing up map $Y \lrarrow \Spec S$. 
For negative $p$, then, the vanishing of 
 $[H^d_{\fN}(S^{\natural})]_p$ follows from the vanishing of 
$H^{d-1}_E(\cO_Y(p))$. By duality, this is equivalent to the
vanishing of $H^1(Y, \omega_Y(-p))$. This vanishing in turn 
is a consequence of
 Lemma \ref{GR} below.

Also, the middle vertical arrow in the commutative diagram above
 is an isomorphism.
Indeed, because $S_{p+q}=S_pS_q$ for all $p,q\gg 0$, we have  
$(\fN^+)^n(S^{\natural}/S^{\dagger})=0$ for $n\gg 0$, whence
 $H_{\fN^+}^i(S^{\natural}/S^{\dagger})=0$ for all $i>0$. 
A look at the long exact cohomology sequence 
corresponding to the exact sequence 
\begin{equation*}0\lrarrow S^{\dagger}
 \lrarrow S^{\natural} \lrarrow S^{\natural}/S^{\dagger}\lrarrow 0
\end{equation*}
then shows that $H_{\fN^+}^i(S^{\dagger})=H_{\fN^+}^i(S^{\natural})$ 
for all $i >1$.

We are now in a position to
 use the snake lemma to compute the cokernel of $\psi$ in
degree $p <0$. For $p < 0$, we get 
 an isomorphism
\begin{equation*}
[\Ker\psi]_{p,q}\cong [\Coker \varphi]_{p,q}.
\end{equation*}
Thus we need to understand the map $\varphi$ in degree $p$ negative.
 
To understand $\varphi$, we first look at  
 the local cohomology modules of
$S^{\dagger}$ and $S^{\natural}$ with supports in $\fN'$ and $\fN''$.
Note that the homogeneous maximal ideal $\fm$ of $S$ 
extends to
$\fN'$ under
 the homomorphism $S\lrarrow S^{\dagger}$ obtained by considering 
  $S$ as the subring of  the bigraded ring $S^{\dagger}$  
concentrated in degrees $(p,0)$, for  $p\ge 0$.
Thus 
   \begin{equation*}H_{\fN'}^d(S^{\natural})=H_{\fm}^d
(S^{\natural})=\bigoplus_{q\ge 0}H_{\fm}^d(S^{\natural}_{\bullet,q})
= \bigoplus_{q\geq0}
H^d_{\fm}(S_{\geq q}),
\end{equation*} where $S^{\natural}_{\bullet,q}$ denotes the graded 
$S$-module $\bigoplus_{p\ge 0}S_{p+q} = S_{\geq q}$.
 Now using
 the long exact sequence of cohomology corresponding to the exact
sequence 
\begin{equation*}
0\lrarrow S^{\natural}_{\bullet,q} \lrarrow S(q)
\lrarrow \bigoplus_{p=-q}^{-1}S_{p+q}\lrarrow 0, 
\end{equation*} we get
 \begin{equation*}
H_{\fN'}^d(S^{\natural})=\bigoplus_{q\ge 0}H_{\fm}^d(S(q)),\end{equation*} 
and so by symmetry also 
\begin{equation*}
H_{\fN''}^d(S^{\natural})=\bigoplus_{p\ge 0}H_{\fm}^d(S(p)).
\end{equation*} 
It is now clear that 
$ H_{\fN'}^d(S^{\natural})$ vanishes in negative $q$-degrees
and that $ H_{\fN''}^d(S^{\natural})$ vanishes in negative $p$-degrees.

So the cokernel of $\varphi$ in degrees  $p < 0$ is isomorphic
to the cokernel of the natural map 
$$
H^d_{\fN'}(S^{\dagger}) \lrarrow H^d_{\fN'}(S^{\natural})
$$
induced by the inclusion of $S^{\dagger}$ in $S^{\natural}$.
Since $H^{d+1}_{\fN'}(S^{\natural}) = H^{d+1}_{\fm}(S^{\natural})$ is zero, 
this cokernel is precisely
$$
H^d_{\fN'}(S^{\natural}/S^{\dagger})
 \cong \bigoplus_{q\geq 0} H^d_{\fm}(S_{\geq q} /S_qS).
$$
Hence, in degree $(p, q)$, we have 
$$
[\Coker \varphi]_{p,q} = [H^d_{\fm}(S_{\geq q}/S_qS)]_{p}.
$$
Now if $S_q \neq 0$, then $S_{\geq q}/S_qS$ is an $S$-module of
dimension strictly less than $d$, and so 
$[H^d_{\fm}(S_{\geq q}/S_qS)]_{p} = 0$. On the other hand, if
$S_q =0$, then 
$$
[H^d_{\fm}(S_{\geq q}/S_qS)]_{p} = [H^d_{\fm}(S_{\geq q})]_p \cong
[H^d_{\fm}(S(q))]_p =
 [H^d_{\fm}(S)]_{p+q}.
$$
Summarizing, we have for negative $p$
\begin{equation*}
[\Coker \varphi]_{p,q}
=
\begin{cases}
{[H_{\fm}^d(S)]_{p+q}}&\text{if $q\in G$;}\\
0&\text{otherwise.}
\end{cases}
\end{equation*}
It thus follows that for positive $p$ 
 there is an exact sequence 
\begin{equation*}
0\lrarrow [\omega_{S^{\natural}}]_p
\lrarrow [\omega_{S^{\dagger}}]_p \lrarrow Q_p\lrarrow 0\end{equation*}
where 
\begin{equation*}
Q_{p,q}=
\begin{cases}
{[\omega_S]_{p+q}}&\text{if $-q\in G$;}\\
0&\text{otherwise.}
\end{cases}
\end{equation*} 
Finally, 
Lemma~\ref{naturalcanonicalmodulelem} below
ensures that the $p$-th graded piece of
 $\omega_{S^{\natural}}$ is $[\omega_S]_{\geq p+1},$
 so this formula for the cokernel $Q_p$ immediately implies the desired
 formula for the positively graded pieces of 
$\omega_{S^{\dagger}}$.
The proof of 
 Theorem \ref{canonicalmodulethm}
 in the case of positive values of $p$ is complete.
\bigskip

To prove Theorem \ref{canonicalmodulethm} for non-positive $p$ 
under the additional hypothesis that $S$ has negative $a$-invariant,
 we use Lemma \ref{negvan} below.
This guarantees that  the module $H^d_{\fN}(S^{\natural})$
vanishes also for non-negative  $p$. We can 
 then use the snake lemma as above to compute the cokernel
of $\varphi$ also for $p\geq 0$. Doing so yields the desired 
conclusion in a similar way.
\end{proof}

\begin{lem}\label{GR}
Let  $S$ be an equidimensional
 section ring of dimension $d \geq 2$ such that 
$\Proj S$ is Cohen-Macaulay. Assume that 
$H^{d-1}_{\fm}(S)$ vanishes in degrees less than $-p$.
Let $Y $ denote $\Proj S^{\natural}$ with its natural polarization.
 Then
$H^{1}(Y, \omega_Y(p)) = 0$ for all $p \geq 0$.
\end{lem}

\begin{proof}
First note that in the important special case where $\Proj S$ has
rational singularities and is defined over a field of characteristic zero,
the result follows immediately from the Grauert-Riemenschneider
vanishing theorem. Indeed, in this case $Y = \Proj S^{\natural}$ has 
rational singularities. Let  $Z \overset{\nu}\to Y$ be a 
resolution. By the Grauert-Riemenschneider vanishing theorem,
the direct image sheaves $R^j\nu_*(\omega_Z\otimes \nu^*\cO_Y(p))
 = R^j\nu_* \omega_Z\otimes \cO_Y(p)$ vanish for $j \geq 1$.
Thus the corresponding spectral sequence degenerates, and
we have an isomorphism
$$
H^i(Z, \omega_Z \otimes \nu^*\cO_Y(p)) \cong 
H^i(Y, \nu_*\omega_Z \otimes \cO_Y(p))
$$
for all $i \geq 0$. Now because $Y$ has rational
singularities, $\nu_*\omega_Z = \omega_Y$ and so
the required vanishing follows from 
standard 
 vanishing theorem applied to the resolution $Z \la \Spec S$
(see \cite{cut}).

We now assume that $\Proj S$ is Cohen-Macaulay and $H^{d-1}_{\fm}(S)$
vanishes in  degrees less than $-p$. We interpret $Y$ as the
 total space bundle over $X$ for the line bundle $\cL^{-1} = \cO_X(-1)$ 
defining the section ring $S$.  The map $f: Y \la X$ is affine and
smooth, with $\omega_Y = f^*\omega_X \otimes f^*\cO_X(1)$.
Computing as in the proof of 
  \cite[Proposition 6.2.1]{HS}, we see that
$$
H^1(Y, \omega_Y(p)) = H^1(Y, f^*(\omega_X(p+1))) = 
H^1(X, f_*\cO_Y \otimes \omega_X(p+1)) = \bigoplus_{n\geq p+1}
H^1(X, \omega_X(n)).
$$
By duality, 
this module vanishes by our assumption that 
 $H^{d-1}_{\fm}(S)$ vanishes in degrees less than $-p$.
\end{proof}

\begin{lem}
\label{naturalcanonicalmodulelem}
 If $S$ is a section ring of an equidimensional projective variety,
 then  for positive $p$,
\begin{equation*}
[\omega_{S^{\natural}}]_p= [\omega_S]_{\ge p+1}.
\end{equation*}
\end{lem}

\begin{proof} It has been proven 
in~\cite[Proposition 6.2.1]{HS} that if $Y=\Proj S^{\natural}$,
then $\Gamma(Y,\omega_Y(p))=[\omega_S]_{\ge p+1}$ for all $p\in \bbZ$.
 We now use
the general fact that we can identify $[\omega_{S^{\natural}}]_p$ 
with $\Gamma(Y,\omega_Y(p))$
for all $p\ge 1$
(see, for example, \cite[Paragraph 2.6.2]{HS}).
\end{proof}

\begin{lem}\label{negvan}
With notation as above, 
suppose that the section ring $S$ has negative $a$-invariant.
Then $H^d_{\fN}(S^{\natural})$ vanishes
 in non-negative  $p$ (and hence $q$) degrees.
\end{lem}

\begin{proof}
Recall that 
for any filtration 
$\left\{I_i\right\}_{i\in \bbN} $ of a ring $S$ we have two fundamental 
short exact sequences
 $$
0 \lrarrow R^+ \lrarrow R \lrarrow S \lrarrow 0 
\quad{\text {and}}
\quad
 0 \lrarrow R^+(1) \lrarrow R \lrarrow G \lrarrow 0,
$$
where $R$ is the Rees ring $\bigoplus_{i \in \bbN} I_i$, $R^+$ is its irrelevant
ideal
$\bigoplus_{i >0} I_i$ , and $G$ is the associated 
graded ring $\bigoplus I_{i}/I_{i+1}$.
In particular, applying this to the natural filtration of $S$,
we get short exact sequences
 $$
0 \lrarrow {(S^{\natural})}^+ \lrarrow S^{\natural} \lrarrow S \lrarrow 0 
\quad {\text {and }}
\quad
0 \lrarrow {(S^{\natural})}^+(1) \lrarrow S^{\natural} \lrarrow S \lrarrow 0, 
$$
where in the first sequence, the module $S$
 is viewed as concentrated in degree zero, while in the second sequence
it has its usual grading. 

The long exact sequence arising from the
first sequence yields an isomorphism
$H^d_{\fN}({(S^{\natural})}^+) \cong 
H^d_{\fN}(S^{\natural})$ in all non-zero $p$ degrees.
The second sequence gives an 
exact sequence
$$
H^d_{\fN}({(S^{\natural})}^+(1)) \lrarrow
H^d_{\fN}(S^{\natural}) \lrarrow H^d_{\fm}(S).
$$
Therefore, 
if $p \geq 0$, our hypothesis gives a surjection
$$
[H^d_{\fN}({(S^{\natural})}^+(1))]_p \lrarrow
[H^d_{\fN}(S^{\natural})]_{p}.
$$
Since 
$[H^d_{\fN}({(S^{\natural})}^+(1))]_{p} \cong 
[H^d_{\fN}({(S^{\natural})}^+)]_{p+1} \cong 
[H^d_{\fN}(S^{\natural})]_{p+1}
$ for all non-negative $p$, we then 
 get  surjections
$$
[H^d_{\fN}(S^{\natural})]_{p+1} \lrarrow
[H^d_{\fN}(S^{\natural})]_{p}
$$
for all $p\geq 0$. This completes the proof since
evidently 
$H^d_{\fN}(S^{\natural})$ vanishes in large degrees.
\end{proof}

\bigskip
We now have all the ingredients we need to prove Theorem \ref{gradform}. 
\begin{proof}[Proof of Theorem \ref{gradform}] 
For each integer $p$,  
set $\Omega'_p
=\Gamma(Y',{I'}^p\omega_{Y'})$,  where 
$I'$ is the ideal of $S$ generated by all elements of degree $N$
and $Y'=\Proj S[I't]$.
We first claim that for $p \geq 1$,
\begin{equation*}
\Omega'_p = 
\bigoplus_{i\in \bbZ; S_i = 0}[\omega_S]_{pN-i}.
\end{equation*}

Note here that $i$ ranges through all integers such that $S_i$ is zero, including
all negative integers. The grading is not given by $i$, but by the inherited
grading from $\omega_S$.

To see why the claim is true, 
for each $p \geq 1$,
 we can identify  $\Omega'_p$ with the $p$-th homogeneous
component of the graded canonical module of $S[I't]$. 
Since $S[I't]=(S^{\dagger})^{(N)}$, we also have  
$\omega_{S[I't]}=(\omega_{S^{\dagger}})^{(N)}$.  
Now by 
 Theorem~\ref{canonicalmodulethm}, we have  
\begin{equation}
\label{Omegaprimeformula}
\Omega_p' = 
[\omega_{S[I't]}]_p=[\omega_{S^{\dagger}}]_{Np}=
\bigoplus_{i \in \bbN; S_i = 0}[\omega_S]_{Np-i}.
\end{equation}  
(The required vanishing needed for Theorem \ref{canonicalmodulethm}
holds by taking  $N$ to be very large.)

We first show the inclusion $\grcore(I\omega_S) \subset \Omega_d'$.
By Lemma \ref{redI'}, $\grcore(I\omega_S) = \core(I'\omega_S)$.
To compute the latter, we appeal to Theorem \ref{mainalg},
which requires that the inclusion
\begin{equation}
\label{eq3}
\Omega_{d-1}'\cap J\omega_S = J\Omega_{d-2}' 
\end{equation} holds
for every homogeneous minimal reduction of $I$.
To verify (\ref{eq3}),
take any 
homogeneous element $w$ in $\Omega'_{d-1} \cap J\omega_S$, say of degree
$n$. 
Since $J$ is necessarily generated by
elements of degree $N$, we can write 
$w=s_1w_1+\ldots +s_r w_r$ where $s_1,\ldots,s_r$ 
are homogeneous elements of degree $N$ in $J$
and $w_1,\ldots,w_r \in [\omega_S]_{n-N}$. 
We want to show that the $w_i$ are in $\Omega_{d-2}'$.
For this we use 
 (\ref{Omegaprimeformula}). 
If
$n\ge (d-1)N+1$, then
 $[\Omega'_{d-2}]_{n-N}= [\omega_S]_{n-N}$ and we are done. On the other hand, 
when $n<(d-1)N+1$, there is nothing to prove
 unless $n=(d-1)N-i$ for some $i$ such that $S_i = 0.$  But in this
case $[\Omega'_{d-2}]_{n-N}=[\omega_S]_{n-N}$, too.
Thus $\grcore(I\omega_S) \subset \Omega_d'$ as needed.

\medskip
The reverse inclusion is an immediate application of
Lemma \ref{van-alg} below (since there are only finitely many positive $i$
 such that $S_i = 0$). 
Note that this lemma implies that when $S_i = 0$ or when $i $ is negative,
we have $[\omega_S]_{Nd - i} \subset J\omega_S$ for every homogeneous
reduction $J$ of $I$. 
\end{proof}

\begin{lem} 
\label{van-alg}
Let $\cL$ be an ample invertible sheaf
on a  projective variety $X$ 
of dimension $d-1$. Fix any 
$N$ so that $\cL^N$ is globally generated and
 so that
 $H^j(X, \cL^{\otimes N(d-1-j)-i}\otimes \omega_X)
 = 0$ for  $j$ in the range $0 < j < d-1$ and
that  $H^j(X, \cL^{\otimes N(d-j)-i}\otimes \omega_X) = 0$ for positive $j$.
 Then $H^0(X, \cL^{\otimes i})=0$ if and only if 
for any set 
$x_1,\ldots,x_d$ of $d$ global generators for $\cL^N$,
the natural inclusion 
\begin{equation*}
\sum_{j=1}^d x_jH^0(X, \omega_X\otimes \cL^{\otimes[N(d-1) - i]})
\subset H^0(X,\omega_X \otimes \cL^{\otimes (Nd-i)})
\end{equation*}
is an equality.
Similarly,  for any $N \gg i$, 
\begin{equation*}
\sum_{j=1}^d x_jH^0(X,  \cL^{\otimes[N(d-1) + i]})
\subset H^0(X, \cL^{\otimes (Nd +  i)})
\end{equation*}
is an equality if and only if $H^{d-1}(X, \cL^{\otimes i}) = 0$. 
\end{lem}

\begin{proof} The first statement  is essentially  ~\cite[Lemma 6.1.2]{HS}.
The second is proved similarly:

 Consider the Koszul
complex determined by the $x_i$ on $X$, where $\cO_X(N)$ denotes $\cL^N$:
$$
0 \la \cO_X(-Nd) \la \dots \la \cO_X(-N)^{\oplus d} 
\overset{(x_1, \dots, x_d)}\lrarrow \cO_X \la 0.
$$
Since the $x_i$ do not all simultaneously vanish on $X$, this
complex is exact. Now tensor with the invertible sheaf $\cO_X(Nd + i)$
to get an exact sequence
$$
0 \la \cO_X(i) \la \dots \la \cO_X(N(d-1) + i)^{\oplus d} 
\overset{(x_1, \dots, x_d)}\to \cO_X(Nd + i) \la 0.
$$
 Using the same argument as in \cite[Lemma 6.1.2]{HS},
we see that the map 
$$
\Gamma(X, \cO_X(N(d-1) + i)^{\oplus d}
 \overset{(x_1, \dots, x_d)}\lrarrow \Gamma(X, \cO_X(Nd + i)) 
$$
is surjective if and only if
$
H^{d-1}(X, \cO_X(i)) 
$
vanishes.
 (Here, we are using 
 that for large enough $N$, 
 $H^j(X, \cO_X(N(d-j) +i )$ vanishes
for $0 < j <d$ and 
$H^{j-1}(\cO_X(N(d-j) + i) $ vanishes for $1< j < d.$)
\end{proof}

\begin{cor}\label{cor1}
Let $S$ be a Cohen-Macaulay section ring, of dimension $d \geq 2$ and
$a$-invariant $a$.  Let $I = S_{\geq N}$ for 
$N \gg 0$. Then 
$$
 S_{\geq Nd + a + 1} \subset 
\grcore(I) \subset
 \bigoplus_{i \in \bbZ; S_i = 0}S_{Nd +a - i}. 
$$
Furthermore, 
the second inclusion is an equality if $S$ is balanced in the
sense that the entire socle of $H^d_{\fm}(S)$ lives in degree $a$.
\end{cor}

\begin{proof}
To prove 
$
\grcore(I) \subset  \bigoplus_{i \in \bbZ; S_i = 0}
S_{Nd +a - i}$, take any $z \in \grcore(I)$. Without loss of 
generality, we may assume that $z$ is homogeneous, say of degree $b$. 
Take any $w \in \omega_S$ of degree precisely $-a$, the minimal possible.
Then $zw \in \grcore(I\omega_S)$ is an element of degree $b-a$.
 Thus by Theorem
\ref{gradform},  $b-a = Nd - i$ for some $i$ such that 
$S_i =0$.  This degree restriction on $z$ forces $z $ in 
 $S_{Nd +a - i}$. Thus the one inclusion is proved.

The  inclusion $S_{\geq dN + a +1} \subset \grcore(I)$ follows 
 immediately from the second part
of Lemma \ref{van-alg}. Fix any homogeneous system of parameters 
$x_1, \dots, x_d$ of degree $N$. To see that 
$S_{Nd+i} \subset (x_1, \dots, x_d)$, we need the vanishing of 
  $H^{d-1}(X, \cO_X(i)) = [H^{d}_{\fm}(S)]_{i},$
where $X = \Proj S$.
For $i \geq a+1$, this 
follows simply from the definition of the $a$-invariant
as the largest $i$ such that $[H^{d}_{\fm}(S)]_{i}$ is non-zero.
When the socle of $H^d_{\fm}(S)$ is concentrated in degree $a$,
Lemma \ref{van-alg} also implies that $S_{Nd + a  - i}
 \subset (x_1, \dots, x_d)S$  whenever $S_i = 0$. Indeed, 
  we need to check only that $ [H^{d}_{\fm}(S)]_{a-i}$ is zero.
But if  $\eta$ is a non-zero element of $ [H^{d}_{\fm}(S)]$,
 then $\eta$ must have a
non-zero multiple in the socle $ [H^{d}_{\fm}(S)]_{a}$.
But if $\eta$ has degree $a-i$, 
this is impossible whenever  $S_i = 0$.
\end{proof}

\medskip

Putting together our Theorem \ref{mainalg} with the results
of \cite{HS}, we arrive at the 
following corollary.
\begin{cor}\label{KC}
Let $\cL$ be an ample invertible sheaf on a rationally singular
variety $X$. Then $\cL$ has a non-zero global section if and
only if
$$
\grcore(I\omega_S) = \adj(I^d\omega_S)
$$
where $S$ is the section ring of the pair $(X, \cL)$,
$d$ is its dimension, I is the ideal of $S$ generated by all 
elements of degrees at least $N$, for $N \gg 0$, 
and $\adj(I^d\omega_S)$ is the adjoint module of $I^d$ as defined
in \cite[Remark 6.0.7]{HS}.
\end{cor}

Thus, using \cite{HS}, Kawamata's Conjecture \cite{kaw1}
is equivalent to the following conjecture.

\begin{conj}\label{KC1}
Let $S$ be the section ring of an ample Cartier divisor $D$ on a 
rationally singular $\bbQ$-Gorenstein variety $X$. Assume
that the $\bbQ$-Cartier divisor $D - K_X$ is ample.
Then 
$$
\grcore(I\omega_S) = \adj(I^d\omega_S),
$$
where $I = S_{\geq N}$, for $N\gg 0$. 
\end{conj}

An interesting special case of Kawamata's Conjecture is
 when $X$ is Fano and  $D = -K_X$. Then the
conjecture collapses to the
following.

\begin{conj}
Let $S$ be a  section ring 
with 
 Gorenstein rational singularities.
Then 
$$
\grcore(I) = \adj(I^d),
$$
where $I = S_{\geq N}$ 
for $N \gg 0$. Here $\adj(I^d)$ denotes the adjoint (or multiplier) ideal
(as defined, for example, in  \cite{PAG}; see also \cite{L3}).
\end{conj}

\bigskip
Finally,
for completeness, we also include the following result 
valid for graded rings of dimension one. The argument is different.

\begin{thm} Let $S$ be a one 
dimensional $\bbN$-graded Cohen-Macaulay ring defined over a field
 containing the rational 
numbers. Suppose that $S$ has homogeneous non-zero-divisors
 of degree $N$ for all  $N\gg 0$. Let $I=S_{\ge N}$ be the 
ideal generated by all elements of degree at least $N$ 
where $N\gg 0$. Then $\grcore(I\omega_S)=\core(I\omega_S)$ if and 
only if there exists a non-zero divisor  of degree one in $S$.
Furthermore, in this case, $\core(I) = \grcore(I)$.  \end{thm}

\begin{proof}
We know by 
\cite[Corollary 5.2.1
and Lemmas 3.1.4 and 3.4.5]{HS}
 that for any homogeneous ideal $I$,
$\core(I\omega_S)=J^{r+1}\omega_S:_{\omega_S} I^r$ for $r\gg 0$,
where $J$ is any reduction of $I$, and $r \gg 0$.

In particular, if $J$ is an ideal of $S$ 
generated by a non-zero divisor  $x$ of degree $N$,
then 
\begin{equation*}
\core(I\omega_S) = J^{r+1}\omega_S:_{\omega_S} I^r \qquad{\text{and}}
\qquad
\grcore(I\omega_S) = J^{r+1}\omega_S:_{\omega_S}{I'}^r,
\end{equation*}
where $I = S_{\geq N}$ is the ideal generated by all elements
of degree at least $N$ and $I' = S_NS$ is the ideal generated
by the elements of degree precisely $N$. 

Now observe that, as submodules of $\omega_S $ tensored with the
full ring of quotients
 of $S$,  
\begin{equation*}
J^{r+1}\omega_S:_{\omega_S} I^r=x^{r+1}\Hom_S(I^r,\omega_S)
\quad {\text{ and }} \quad
J^{r+1}\omega_S:_{\omega_S} {I'}^r=x^{r+1} \Hom_S({I'}^r,\omega_S).
\end{equation*}
 So $\grcore(I\omega_S)=\core(I\omega_S)$ if 
and only if $\Hom_S(I^r,\omega_S)=\Hom_S({I'}^r,\omega_S)$.
 But as $S$ is  Cohen-Macaulay,
this is  equivalent to $I^r={I'}^r$.

Suppose now that $\grcore(I\omega_S)=\core(I\omega_S)$. 
Then $I^r={I'}^r$, and so  $I^r$ is generated in degrees precisely
$Nr$. In particular, 
$S_{rN+1}=S_1 S_{rN}$. This means that 
 $S_{rN+1}$ is contained in the ideal generated by  $S_1$.
Now 
 as
there are non-zero-divisors
 of degree $rN+1$,  there must be a non-zero-divisor of degree one.
Otherwise, $S_{Nr+1}$ is contained in a minimal prime 
and consists only of zero-divisors.

Conversely, if $S$ admits a non-zero-divisor $y$ of degree one, 
then there exists a $q$ such that $S_{\geq n} = y ^{n-q}S_{\geq q}$ for
all $n \geq q$. In particular, for each $p \geq 0$ and $N \geq q$, 
$$
[I^r]_{Nr+p} \subset  S_{Nr + p} \subset y^{Nr + p - q}S_q \subset
((y^N)^{r-1} y^{N-q}S_q)S \subset (I')^r.
$$
Since $I^r = \oplus_{p \geq 0} [I^r]_{Nr + p}$, 
we conclude that $I^r = {I'}^r$ and so the core and graded core modules
 are equal. 

For the final statement, assume that $\core(I\omega_S) = \grcore(I\omega_S),$
and let $z \in \grcore(I)$. Then $z \omega_S \in \grcore(I\omega_S) =
 \core(I\omega_S),$ so that for every reduction $J$ of $I$, we have
$z \in (J\omega_S:_S \omega_S).$ As any such  $J$ is generated
 by a non-zero-divisor,
it follows that $z \in J$ for all reductions $J$ of $I$,
 and so $z \in \core(I)$.
\end{proof}

\section{The core of powers of $\fm$ in a standard graded ring}

We now turn to the proof of Theorem \ref{standard} from the introduction,
giving a formula for core in a reduced standard graded ring of
 arbitrary characteristic.

\begin{thm}\label{standard2}
Let $(S, \fm)$
 be a reduced Cohen-Macaulay  graded ring generated in degree one
over an infinite field.
Let $d$ denote the dimension of $S$ and let 
$a$ denote its $a$-invariant.
Then 
$$
\core(\fm^n) = \grcore(\fm^n) =  \fm^{nd +a +1}
$$
for all $n \in \NN$.
\end{thm}

Note that there is {\it no restriction} on the characteristic of the ground
field on Theorem \ref{standard2}.
This is especially interesting in light of examples from  \cite{PU},
which show in general that the core depends on the characteristic.

\begin{proof}
By \cite[Theorem 4.5]{CPU}, 
we know that $\core(\fm^n) = \grcore(\fm^n),$ so we consider only the
graded core.
 Lemma \ref{van-alg} shows that 
$\fm^{nd + a + 1} \subset \grcore(\fm^{n}),$ for all $n$ (not just large $n$). 
(Note that because $S$ has depth at least two, $S_n = H^0(X, \cO_X(n))$
for all $n$.)
It remains to prove the other inclusion.

The point is to prove the analogous statement
inside the canonical module of $S$, namely that 
\begin{equation}\label{4}
\grcore(\fm^n\omega_S) \subset [\omega_S]_{\geq nd +1}.
\end{equation}
This will imply that $\grcore(\fm^n) \subset S_{\geq nd +a +1}$  
as follows. Suppose that
 $z\in \grcore(\fm^n)$ has degree less than $nd+a+1$. 
Choose any $w \in \omega_S$ of degree $-a$. Then $zw \in 
 \grcore(\fm^n\omega_S)$ is an element of degree less than $nd +1$,
 contradicting (\ref{4}).

To see that 
$\grcore(\fm^n\omega_S) \subset [\omega_S]_{\geq nd +1},$
we use Lemma \ref{key} below applied to the 
ideal $I = \fm$.  The hypothesis is satisfied by easy degree reasons: since
$\Omega_p = [\omega_S]_{\geq p+1},$
and $J$ is generated by elements of degree $n$,
it follows immediately   that $ \Omega_{nd-n} \cap  J\omega_S = 
J \Omega_{nd-2n}$.
Note that the associated
graded ring $G$ is simply again $S$.
So to show that 
$$\grcore(\fm^n) \subset \Omega_{nd} = [\omega_S]_{\geq nd+1}$$
it suffices to show the vanishing of  
$$
\left[\bigcap(x_1^*, \dots, x_d^*)\omega_G\right]_{nd} = 
\left[\bigcap(x_1, \dots, x_d)\omega_S\right]_{nd},
$$
where the intersection is taken over all homogeneous systems
of parameters of degree $n$.
Now letting $T = S^{(n)}$ denote the $n$-th Veronese
subring of $S$, we have that 
$$
\left[\bigcap(x_1, \dots, x_d)\omega_S\right]_{nd} = 
\left[\bigcap(x_1, \dots, x_d)\omega_{T} \right]_{d}.
$$ 
Thus we need to prove the vanishing of $\cap J\omega_T$ in degree $d$ as
$J$ ranges over all homogeneous systems of parameters
of $T$ of degree one. But since 
$T$ 
 is a reduced 
Cohen-Macaulay ring of dimension $d$ generated in degree one, 
the vanishing of this intersection is an immediate consequence of 
Lemma  \ref{reduced} below. 
\end{proof}

\begin{lem}\label{key}
Let $I \subset A$ be an $\fm$-primary ideal in a local or graded ring $A$
of   dimension $d$ at least two,
and let $G$ denote the associated graded ring of $A$ with respect to $I$.
Fix any $n \in \NN$,
and let $\cR$ be any non-empty set of reductions of $I^n$.
Assume that $\Omega_{nd-n} \cap J\omega_A = J \Omega_{nd-2n}$
for all $J \in \cR$. 
Then 
$$
\bigcap_{\{x_1, \dots, x_d\} \in \cR} (x_1, \dots, x_d)\omega_A \subset 
\Omega_{dn} = \Gamma(\Proj A[It], I^{nd}\omega_{\Proj A[It]})  
$$
if 
\begin{equation}
\left[\bigcap_{\{x_1, \dots, x_d\} \in \cR} (x_1^*, \dots, x_d^*)\omega_G
\right]_{nd} = 0,
\end{equation}
where $x_i^*$ denotes the image of $x_i$ in $I^n/I^{n+1}$,
considered as an element of $G$. The converse holds in $G$ is Cohen-Macaulay.
\end{lem}

\begin{proof} This is a slight variant of the Key Lemma,
 \cite[Lemma 3.3.1]{HS}. The proof is essentially the same so we omit it.
\end{proof}

\begin{lem}\label{reduced}
Let $T$ be  a reduced Cohen-Macaulay graded ring  of positive
 dimension 
$d$ finitely generated by its
degree one elements over an infinite
field $k$. Then $
\bigcap (x_1, \dots, x_d)\omega_T 
$ vanishes in degree $d$, where the intersection is taken 
over all homogeneous systems of parameters of degree one.
\end{lem}

\begin{proof}
By induction on $d$
(see \cite[Proof of Theorem 3.4.1]{HS}),
we reduce immediately to the case where $T$ has dimension one.
Note that in doing this induction, we can kill $d-1$ generic one-forms,
so that the resulting one dimensional ring can be assumed to be reduced.  
The scheme $\Proj T$ is a finite set of  reduced
 points in some projective
space $\bbP^n = \bbP(T_1^*)$. Say there are $e$ points.

Now, by Lemma \ref{equiv1} below, it suffices to prove the dual statement, 
namely that  the module $H^1_{\fm_T}(T)$ is generated in degree $-1$ 
as a $k$-vector space by elements of the form
$\left[\frac{1}{x}\right]$, where $x$ ranges through non-zero divisors 
in 
 $T$
 of degree one.

Now let $y$ be 
 any non-zero-divisor of degree $e$. Then thinking of $H^1_{\fm_T}(T)$
as the direct limit $T/y \lrarrow T/y^2 \lrarrow \dots$ whose maps are
given by multiplication by $y$, it is easy to check that 
$$
[H^1_{\fm_T}(T)]_{-1} = \left[\frac{T_{e-1}}{y}\right].
$$
We now claim that we can find non-zero-divisors $x_1, \dots, x_e$ of degree
one such that 
$$
\left[\frac{1}{x_1}\right], \dots, \left[\frac{1}{x_e}\right]
$$
generate $[H^1_{\fm_T}(T)]_{-1}$ as a $k$-vector space.
This is equivalent to asking that the elements
$$
\left[\frac{x_1\dots \hat x_i \dots x_e}{x_1 \dots x_e}\right] \qquad {\text{ 
for }} i = 1, \dots e
$$
generate $[H^1_{\fm_T}(T)]_{-1}$,
which is in turn equivalent to 
asking that the forms 
$$
\ell_i = x_1\dots \hat x_i \dots x_e \qquad {\text{ for }} i = 1, \dots, e
$$
span $T_{e-1}$. (Here, the notation $\hat x_i$ 
means that $x_i$ has been omitted.)

To find such $x_i$, just let $x_i$ be any linear form on $\bbP^n$
vanishing at the point $P_i$ but not at any other point of $\Proj T$.
Then $\ell_i$ vanishes at all the points except $P_i$
and so the  $\ell_i$ are linearly independent elements of $T_{e-1}$.
On the other hand, since the dimension of $T_{e-1}$ is $e$, 
these elements $\ell_i$ span $T_{e-1}$.
The proof is complete. 
\end{proof}

\begin{lem}\label{equiv1}
Let $G$ be any $\bbN$-graded ring generated in degree one as an algebra
over an Artinian local ring $G_0$.  Let $d$ denote the dimension of $G$.
Let $\cR$ be any set of ideals generated by homogeneous systems of 
parameters for $G$ of degree one. Then
$$
\bigcap_{J\in \cR} J\omega_G = 0$$
if $
[H^d_{\fm_G}(G)]_{-d}\,\,
$
is generated  over $G_0$  by elements of the form 
$\left[\frac{1}{x_1\dots x_d}\right]$, for $(x_1, \dots, x_d) \in \cR$.
The converse holds if $G$ is Cohen-Macaulay.
\end{lem}

\begin{proof}
We use the natural duality pairing
\begin{equation*}
[H_{\fM_G}^d(G)]_{-d}\otimes_A [\omega_G]_d \lrarrow [E_G(k)]_0  =
 E_{G_0}(k) = \omega_{G_0}, 
\end{equation*}
where $k$ is the residue field of $G $  (and $G_0$), 
sending
$$
\gh \otimes \omega \mapsto 
\omega(\gh).
$$
More explicitly, making the identification, $E_G(k)=
H_{\fM}^d(\omega_G)$, we have
\begin{equation*}
\Big[\frac{g}{(x_1^*\cdots x_d^*)^t}\Big]\otimes \omega
\mapsto \Big[\frac{g\omega}{(x_1^*\cdots x_d^*)^t}\Big].
\end{equation*}

Now, let us first suppose that $[H^d_{\fM_G}(G)]_{-d}$
is generated as an $G_0$-module
by the special elements $\Big[\frac{1}{x_1^*\cdots x_d^*}\Big]$.
Take any 
$w \in \bigcap(x^*_1, \dots, x_d^*)\omega_G$ of degree $d$.
Then $\langle \eta, w \rangle = 0$ for any element $\eta$ 
of the form 
$\Big[\frac{1}{x_1^*\cdots x_d^*}\Big]$. So since these generate
$[H^d_{\fM_G}(G)]_{-d}$, we see that 
$\langle \eta, w \rangle = 0$ for all $\eta \in [H^d_{\fM_G}(G)]_{-d}$.
Since the pairing is perfect this implies that in fact $w$ is zero,
as needed.

Conversely, suppose that $ \bigcap(x^*_1, \dots, x_d^*)\omega_G$
contains no non-trivial elements of degree $d$.
Let  $H$ be the $G_0$-submodule of $H^d_{\fM_G}(G)$ generated by 
the special elements 
 $\Big[\frac{1}{x_1^*\cdots x_d^*}\Big]$.
Consider
the exact sequence \begin{equation*}0\lrarrow H \lrarrow
[H_{\fM}^d(G)]_{-d}\lrarrow C\lrarrow 0.\end{equation*}
Taking the Matlis dual over the zero dimensional ring $G_0$,
we obtain
\begin{equation*}\Hom_{G_0}(C, \omega_{G_0})=\Ker\Big([\omega_G]_d
\lrarrow \Hom_{G_0}(H, \omega_{G_0})\Big).
\end{equation*}
We now compute the kernel of the map 
$[\omega_G]_d
\lrarrow \Hom_{G_0}(H,\omega_{G_0}).$
An element $w$ is sent to the zero map $H \lrarrow \omega_{G_0}$
if and only if for all $\eta \in H$, the pairing 
$\langle \eta, w \rangle = 0$. But this occurs if and 
only if for all $(x_1, \dots, x_d) \in \cR$, 
the element 
$\Big[\frac{w}{x_1^*\cdots x_d^*}\Big]$ is zero.
Since $\omega_G$ is Cohen-Macaulay, this in turn means that $w \in
 (x_1^*, \dots, x_d^*)\omega_G$ for all $(x_1, \dots, x_d) $ in $\cR$. The
proof is complete.
\end{proof}

\section{Questions}

In this final section, we gather together a few commutative algebra
questions whose positive solutions would solve Kawamata's Conjecture,
at least in special cases. To keep things simple, we will focus
on the case where $H^i(X, \cO_X)$ vanishes for $i>0$, which simplifies 
statements considerably. This  has the effect of allowing us to assume
that the section rings $S$ have rational singularities (see  \cite[\S6]{HS}).
All of these questions are interesting if we furthermore assume that
$S$ is Gorenstein with isolated singularity and
 $a$-invariant -1; this corresponds to 
the important special  case where $\cL$ is the anti-canonical divisor on
a smooth Fano variety.

\begin{question}\label{q1}
Suppose that $ S$ is a section ring with  rational singularities and let
$Y' \lrarrow \Spec S$ be the proper birational map
obtained by blowing up the ideal generated by all elements of degree $N$, 
$N \gg 0$.   Let $\nu: Y \lrarrow Y'$
be the normalization map. Then  is $\nu_*\omega_{Y} = \omega_{Y'}$?
Note that the rational singularities hypothesis implies that at least 
there is an equality for their 
 global sections.
\end{question}

An affirmative answer to Question \ref{q1} 
 is equivalent to the
existence of a non-zero element in $S$ of degree one.
Indeed, we know that the core module of $I$ is 
$$
\Omega_d = \Gamma (Y, I^d\omega_Y) = \Gamma (Y, (I')^d\cO_Y \otimes
\omega_Y) =
\Gamma(Y',   (I')^d\cO_{Y'} \otimes \nu_*\omega_Y),
$$
 so if $\nu_*\omega_{Y} = \omega_{Y'}$, then
this is also equal to 
$$
\Gamma(Y',  (I')^d \otimes\omega_{Y'}) = \Omega'_d,
$$ 
which is the graded core module of $I$.

 The map $Y \lrarrow Y'$ is an isomorphism except along a
single irreducible divisor 
$E'$  isomorphic to $\Proj S$.
Furthermore, because $\nu_*\omega_Y$ is $S_2$, it suffices to check the
isomorphism $\nu_*\omega_{Y} \lrarrow \omega_{Y'}$ generically along $E'$.
Note also that the preimage of $E'$ in $Y$ is a single irreducible divisor $E$
also isomorphic to $\Proj S$.

\medskip
 Question \ref{q1} can be rephrased more algebraically
as follows:

\begin{question}\label{q2}
Let $S$ be a section ring with rational singularities.
Let 
 $I'$ be the ideal generated by all elements of degree $N$, for $N \gg 0$,
and 
consider the normalization map $S[I't] \subset S[It]$.
Then is the  naturally induced trace map of graded canonical modules 
\begin{equation}\label{5}
\omega_{S[It]} \subset \omega_{S[I't]}
\end{equation}
an isomorphism? Equivalently,
 is it an isomorphism in sufficiently high degrees, or
after localizing at the height one prime containing $I'S[I't]$?
\end{question}

Note that
 when $S_1 \neq 0$, the map  (\ref{5}) is an isomorphism by Theorem
\ref{canonicalmodulethm}. Thus in this setting, settling Question 1 or 2 
 is essentially equivalent to Kawamata's Conjecture.

\medskip
It is worth pointing out that in Questions 1 and 2, we cannot assume,
nor can we hope to show in general, that $Y'$ satisfies Serre's $S_2$ 
condition:

\begin{prop}
Let $S$ be a section ring of a normal variety
with respect to an ample invertible sheaf $\cL$,
 and assume that $S_1 \neq 0$.
Let $I' $ be the ideal generated by elements of degrees precisely $N$,
for $N \gg 0$.
Then if $\Proj S[I']$ satisfies Serre's $S_2$-condition, then $S_1$ 
generates an $\fm$-primary ideal. In particular,
$\cL$ is  globally generated. 
\end{prop}

\begin{proof}
By Theorem \ref{canonicalmodulethm} and Lemma \ref{naturalcanonicalmodulelem},
we see that if $S_1 \neq 0$, the inclusion $\omega_{S^{\natural}} \subset 
\omega_{S^{\dagger}}$ is an isomorphism in positive degrees. This forces
$\nu_*\omega_Y \cong \omega_{Y'}$. In particular, we have an isomorphism
 $(\nu_*\omega_Y)_{E'} \cong (\omega_{Y'})_{E'}$ along the reduced exceptional
divisor $E'$. This means that the induced trace map for 
the map of one-dimensional local
rings $\cO_{Y',E'} \hookrightarrow \cO_{Y,E}$ is an isomorphism.
Let $Q$ be the cokernel of the injection 
 $\cO_{Y',E'} \hookrightarrow \cO_{Y,E}$; note that $Q$ has
 dimension zero.
 We have an induced exact sequence of
local cohomology modules 
$$
0 \lrarrow Q =  H^0_{\fm}(Q) \lrarrow H^1_{\fm}(\cO_{Y', E'})
\lrarrow H^1_{\fm}(\cO_{Y,E}) \lrarrow 0
$$
where $\fm$ is the  maximal ideal of $\cO_{Y', E'}$.
Dualizing, we see that the dual of 
$Q= H^0_{\fm}(Q)$ is the cokernel of the trace map.
Hence if the trace map in an isomorphism, then $Q$ is zero and
the  map $\cO_{Y'E'} \hookrightarrow \cO_{Y,E}$ is an isomorphism.
This shows that $Y'$ is regular along $E'$, and hence regular everywhere
in codimension one. 

Now if $Y'$ is also $S_2$, then it must be normal, whence $Y = Y'$.
But then since $Y = \Proj S[It]$ and $Y' = \Proj S[I't]$,
it follows that $I^r = {I'}^r$ for $r $ sufficiently large. 
Thus $S_{\geq Nr} = S_{Nr}S$ for large $r$.
In particular,
$S_{Nr +1} = S_{Nr}S_1$, which means that the 
ideal generated by 
$S_{Nr+1}$ is contained in the ideal generated by $S_1$.
Since $S_{Nr+1}$ generates an $\fm$-primary ideal,
it follows that so does $S_1$. The proof is complete.
\end{proof}

This proof shows that the isomorphism $\nu_*\omega_Y \cong \omega_{Y'}$
holds if and only if the induced map of one dimensional
local rings
$\cO_{Y', E'} \hookrightarrow \cO_{Y, E} $ is an isomorphism.
Since this is the normalization map, we see that $\cL$ has a section
if and only of the scheme 
 $Y'$ is regular in codimension one. Thus we can ask:

\begin{question}\label{q4}
If $S$ is a rationally singular section ring and $I'$ is the ideal 
generated by all elements of degrees precisely $N$, for $N \gg 0$,
then is $\Proj S[I't]$ regular in codimension one?
Equivalently, is the local ring $S[I't]_P$ regular where
$P$ is the unique height one prime containing $I'S[I't]$?
\end{question}

Another possible approach to Question 1 is the following.

\begin{question}
Let $S$ be a rationally singular section ring of dimension $d \geq 2$.
Fix $N \gg 0$.
Let $I'$ be the 
ideal generated by elements of degree precisely $N$, and let $L$ 
be the ideal generated by all elements of degrees at least $Nd$.
Consider the multi-Rees ring $R = S[I't][Lu]$, and let $\fM$ denote
its unique homogeneous maximal ideal. Then does
$H^{d+1}_{\fM}(R)$ vanish? Equivalently, does it vanish in $u$-degree zero?
\end{question}

Interestingly, the trace map for the normalization map of  $R$ is an 
isomorphism. So $\omega_R = \omega_{S[It][Lu]}$ and in particular,
$R$ has a Cohen-Macaulay canonical module.

An affirmative solution to Question 4 is equivalent to an affirmative answer to
Question 1. 
To see this,  first note that for $N\gg 0$, we can assume that 
cohomology module $H^{d+1}_{\fM}(R)$ is concentrated in $u$-degree zero.
Now the point is essentially that $Y$ is obtained from $Y'$
 by blowing up the sheaf of
ideals $L\cO_{Y'}$. One specific way to see the equivalence is to note that 
 the $u$-degree
 zero part of $H^{d+1}_{\fM}(R)$ is precisely the dual of the
cokernel of the trace map $\Gamma(\omega_{\Proj R}) \lrarrow
 \omega_{\Spec S[I't]}$, where the scheme $\Proj R$ is
 defined via the $u$-grading on $R$. But now it is easy to check that 
 $\Proj R \cong \SSpec_Y \Sym (\cO_{\Proj{S[Lt]}}(1))$, so that there
is a smooth map $\Proj R \lrarrow Y$ realizing $\Proj R$ as a line bundle 
over the smooth variety $Y = \Proj S[Lt]$.
Thus, one computes directly (as in \cite[Proposition 6.2.1]{HS}) that
 $\Gamma(\omega_{\Proj R})$
is $
\bigoplus_{p \geq 1} \Gamma(Y, I^p\omega_Y).
$ 
So the trace $\Gamma(\omega_{\Proj R}) \lrarrow
 \omega_{\Spec S[I't]}$ is 
$$
\bigoplus_{p \geq 1} \Gamma(Y, I^p\omega_Y) = \bigoplus_{p \geq 1}
 \Omega_p \hookrightarrow
\bigoplus_{p \geq 1} \Gamma(Y', {I'}^p\omega_{Y'})
= \bigoplus_{p\geq 1} \Omega_p'.
$$
But as we know, the trace map $\nu:Y \lrarrow Y'$ is
an isomorphism if and only if we have equalities 
$\Omega_p' = \Omega_p$
for all $p \geq 1$.

\bigskip
We now consider a slightly different direction.
\begin{question}
Let $S$ be a rationally singular section ring of dimension $d \geq 2$.
Set $S^{\dagger} = \bigoplus_{p \in \bbN} S_{p}S$ as in Paragraph \ref{rings}.
Is the $a$-invariant of $S^{\dagger}$ negative?
That is, does $\omega_{S^{\dagger}}$ vanish in non-positive degrees?
\end{question}

The interest in Question 3 stems from the following.
\begin{prop} Let $S$ be an equidimensional
 section ring of dimension $d\ge 2$ with an isolated non-Cohen-Macaulay point.
 Then $S_1\not=0$ if and only if
$a(S^{\dagger})<0$. 
\end{prop}

\begin{proof}
Consider $\omega_{S^{\dagger}}$ as a bigraded module;
its homogeneous components are given by Theorem \ref{canonicalmodulethm}.
In particular, if
 $S_1=0$, then for $p > 0$, 
the degree $(p,-1)$ component of 
$\omega_{S^{\dagger}}$ is 
 $[\omega_S]_{p-1}$, and so in particular is non-zero.
 By $p$-$q$ symmetry, then also the degree 
$(-1, p)$ component of
$\omega_{S^{\dagger}}$ is 
 non-zero. Thus $\omega_{S^{\dagger}}$ does not vanish in $p$-degree $-1$,
and so $a(S^{\dagger}) \geq 1$.

Conversely, if $a(S^{\dagger}) < 0$, then 
$\omega_{S^{\dagger}}$ vanishes in  $p$ degree -1. In particular,
$[\omega_{S^{\dagger}}]_{-1, q}$ is zero and by symmetry then 
$[\omega_{S^{\dagger}}]_{q, -1}$ is zero as well.
From Theorem \ref{canonicalmodulethm}, it follows that 
$S_1$ must be non-zero. The proof is complete.
\end{proof}

\bigskip
Finally, we mention one more group of questions,  
 stemming from our approach to computing core in \cite{HS}.
First a bit a notation:
For any $\fm$-primary ideal in a local (or graded) 
 ring $(A, \fm)$ of dimension $d >0$, let $G = A/I \oplus I/I^2 \oplus \dots$
be the associated graded ring of $I$. For any element $x \in I \setminus I^2$,
let $x^*$ denote the  class of $x$ in $I/I^2$, considered as
a degree one element of $G$. Elements of the local cohomology
module $H^d_{\fM_G}(G)$ can be denoted in the usual way as 
generalized
fractions 
\begin{equation*}\Big[\frac{g}{(x_1^*\cdots x_d^*)^t}\Big]\end{equation*}
where $g\in G$ and $t\in \Bbb N$, and $x_1, \dots, x_d$ is a minimal reduction
of $I$. Note that the degree of such an element is
$\deg{g} - dt$.

\begin{question}\label{q6}
Let $S$ be a rationally singular section ring and let $I$ be the ideal
generated by all elements of degrees at least $N$, for $N \gg 0$.
Let $
G = S/I \oplus I/I^2 \oplus \dots
$
be the associated graded ring of $S$ with respect to $I$.
Then is 
$
[H^d_{\fM_G}(G)]_{-d}
$
generated as an $G_0$-module  by elements of the
form
$$
 \Big[\frac{1}{x_1^*\cdots x_d^*}\Big],
$$
where $x_1, \dots x_d$ ranges through all
 homogeneous system of parameters
of degree $N$?
\end{question}

Lemma \ref{equiv} below  motivates this question. 
Indeed, this lemma implies that if Question \ref{q6}
has an affirmative answer, then 
$$
\grcore(I\omega_S) = \Omega_d = [\omega_S]_{\geq Nd+1},
$$
so by Theorem \ref{gradform}, $S_1$ must be non-zero.
Note that in the setting of Question 6, 
the hypothesis of Lemma \ref{equiv} holds by 
 \cite[Proposition 4.14]{HS}. 
Interestingly,
 the main results of \cite{HS} imply that Question \ref{q6}
does have an affirmative solution if we allow $x_1, \dots, x_d$
to range over {\it all} reductions of $I$, not just homogeneous ones.

\begin{lem}\label{equiv}
Let $I$ be an $\fm$-primary ideal in a local (or graded)
ring $(A, \fm)$ of  dimension $d \geq 2$,
and let $G$ be the corresponding associated graded ring $A/I \oplus I/I^2
\oplus \dots.$ 
Let $\cR$ be any set of minimal reductions of $I$,
and assume that $\Omega_{d-1} \cap J\omega_A = J\Omega_{d-2}$
for all $J \in \cR$.
Then 
$$
 \bigcap_{J \in \cR} J\omega_A \subset \Omega_d 
$$ 
if
$
[H^d_{\fM_G}(G)]_{-d}$
is generated as an $A/I$-module
 by elements of the form 
$$
\Big[\frac{1}{x_1^*\cdots x_d^*}\Big]
$$
where 
$(x_1, \dots, x_d)$ ranges through all the reductions in 
$\cR$.
The converse holds if $G$ is Cohen-Macaulay.
\end{lem}

\begin{proof}
Using Lemma \ref{equiv1}, this is 
 essentially a dual form of the Key Lemma of \cite[Lemma 3.3.1]{HS},
which states that under the same hypothesis,
\begin{equation}\label{!}
\bigcap_{J \in \cR} J\omega_A = \Omega_d 
\quad{\text{ if }}
\quad
\bigcap_{(x_1, \dots, x_d) \in \cR} [(x_1^*, \dots, x_d^*)\omega_G]_{d} = 0,
\end{equation}
and conversely if $G$ is Cohen-Macaulay.
\end{proof}

\bigskip
One approach to Question \ref{q6} would be to find an affirmative answer 
to the next question.

\begin{question}{\label{q5}}
Let $S$ be a rationally singular section ring of dimension $d\geq 2$, and
let $I$ be the ideal generated by elements of degree exactly $N$,
 for $N \gg 0$.
Let $$
F = S/\fm \oplus I/\fm I \oplus I^2/\fm I^2 \oplus \dots
$$
be the fiber ring of $I$ (whose associated projective scheme
defines the scheme-theoretic closed fiber of the map $\Proj S[It] \lrarrow 
\Spec S$), and let $F_{red}$ denote $F$ modulo its nilradical.
Is the map of local cohomology
$$
H^d_{\fm_F}(F) \lrarrow H^d_{\fm_F}(F_{red})
$$
a bijection in degree $-d$?
\end{question}

To see the relationship with Question \ref{q6}, first note that 
$H^d_{\fM_F}(F) \cong H^d_{\fM_G}(G) \otimes_S S/\fm$, 
where $\fm$ is the maximal ideal of $S$. Thus by Nakayama's Lemma,
the required elements generate $[H^d_{\fM_G}(G)]_{-d}$ as an $A/I$-module
if and only if their images modulo $\fm$ generate  
$[H^d_{\fM_F}(F)]_{-d}$ over $k$. So in Question \ref{q6}, we might
as well ask the analogous question for $F$. On the other hand,
it is easy to check that $F_{red}$ (and $G_{red}$) is isomorphic
to  $S^{(N)}$, the $N$-th Veronese subring of the section ring $S$. 
Now because $F_{red}$ is a reduced standard graded Cohen-Macaulay algebra,
Lemma \ref{reduced} ensures that an affirmative answer to Question
\ref{q5} will imply one for Question \ref{q6}.

\newcommand{\noopsort}[1]{} \newcommand{\printfirst}[2]{#1}
  \newcommand{\singleletter}[1]{#1} \newcommand{\switchargs}[2]{#2#1}
\providecommand{\bysame}{\leavevmode\hbox to3em{\hrulefill}\thinspace}
\providecommand{\MR}{\relax\ifhmode\unskip\space\fi MR }
\providecommand{\MRhref}[2]{%
  \href{http://www.ams.org/mathscinet-getitem?mr=#1}{#2}
}
\providecommand{\href}[2]{#2}

\end{document}